\documentclass[10pt,twoside]{article}
\usepackage{Latex-document}
\usepackage{amsmath,amsthm,amsfonts,eucal,epsfig}
\usepackage[matrix,arrow]{xy}

\newcommand{\Q}{\mathbb{Q}}
\newcommand{\C}{\mathbb{C}}

\newcommand{\iso}{\cong}           

\newcommand{\CP}[1]{\C {\mathrm P}^{#1}}

\newtheoremstyle{my}{1.5em}{0.5em}{\em}{}{\sc}{.}{0.5em}{}
\newtheoremstyle{mydef}{1.5em}{0.5em}{}{}{\sc}{.}{0.5em}{}

\theoremstyle{my}
\newtheorem{thm}{Theorem}

\newtheorem{conjecture}[thm]{Conjecture}

\newtheorem{assumption}[thm]{Assumption}

\newcommand{\A}{\mathcal A}
\newcommand{\E}{\mathcal E}
\newcommand{\F}{\mathcal F}
\newcommand{\K}{\mathcal K}
\renewcommand{\O}{\mathcal O}

\begin{document}

\title{Fukaya categories and deformations}
\author{Paul Seidel}
\date{June 15, 2002}
\maketitle

Soon after their first appearance \cite{fukaya93}, Fukaya categories
were brought to the attention of a wider audience through the
homological mirror conjecture \cite{kontsevich94}. Since then Fukaya
and his collaborators have undertaken the vast project of laying down
the foundations, and as a result a fully general definition is
available \cite{fooo,fukaya01b}. The task that symplectic geometers
are now facing is to make these categories into an effective tool,
which in particular means developing more ways of doing computations
in and with them.

For concreteness, the discussion here is limited to projective
varieties which are Calabi-Yau (most of it could be carried out in
much greater generality, in particular the integrability assumption
on the complex structure plays no real role). The first step will be
to remove a hyperplane section from the variety. This makes the
symplectic form exact, which simplifies the pseudo-holomorphic map
theory considerably. Moreover, as far as Fukaya categories are
concerned, the affine piece can be considered as a first
approximation to the projective variety. This is a fairly obvious
idea, even though its proper formulation requires some algebraic
formalism of deformation theory. A basic question is the
finite-dimensionality of the relevant deformation spaces. As
Conjecture \ref{conj:hh} shows, we hope for a favourable answer in
many cases. It remains to be seen whether this is really a viable
strategy for understanding Fukaya categories in interesting examples.

Lack of space and ignorance keeps us from trying to survey related
developments, but we want to give at least a few indications. The
idea of working relative to a divisor is very common in symplectic
geometry; some papers whose viewpoint is close to ours are
\cite{ionel-parker98, li-ruan98, eliashberg-givental-hofer02,
ozsvath-szabo01}. There is also at least one entirely different
approach to Fukaya categories, using Lagrangian fibrations and Morse
theory \cite{fukaya-oh98, kontsevich-soibelman00, fukaya02}. Finally,
the example of the two-torus has been studied extensively
\cite{polishchuk-zaslow98}.

{\em Acknowledgements.} Obviously, the ideas outlined here owe
greatly to Fukaya and Kontsevich. The author is equally indebted to
Auroux, Donaldson, Getzler, Joyce, Khovanov, Smith, and Thomas (an
incomplete list), all of whom have influenced his thinking
considerably. The preparation of this talk at the Institute for
Advanced Study was supported by NSF grant DMS-9729992.

\section{Symplectic cohomology}

We will mostly work in the following setup:

\begin{assumption} \label{as:projective}
$X$ is a smooth complex projective variety with trivial canonical
bundle, and $D$ a smooth hyperplane section in it. We take a suitable
small open neighbourhood $U \supset D$, and consider its complement
$M = X \setminus U$. Both $X$ and $M$ are equipped with the
restriction of the Fubini-Study K{\"a}hler form. Then $M$ is a
compact exact symplectic manifold with contact type boundary,
satisfying $c_1(M) = 0$.
\end{assumption}

Consider a holomorphic map $u: \Sigma \rightarrow X$, where $\Sigma$
is a closed Riemann surface. The symplectic area of $u$ is equal (up
to a constant) to its intersection number with $D$. When counting
such maps in the sense of Gromov-Witten theory, it is convenient to
arrange them in a power series in one variable $t$, where the $t^k$
term encodes the information from curves having intersection number
$k$ with $D$. The $t^0$ term corresponds to constant maps, hence is
sensitive only to the classical topology of $X$. Thus, for instance,
the small quantum cohomology ring $QH^*(X)$ is a deformation of the
ordinary cohomology $H^*(X)$.

As we've seen, there are only constant holomorphic maps from closed
Riemann surfaces to $M = X \setminus D$. But one can get a nontrivial
theory by using punctured surfaces, and deforming the holomorphic map
equation near the punctures through an inhomogeneous term, which
brings the Reeb dynamics on $\partial M$ into play. This can be done
more generally for any exact symplectic manifold with contact type
boundary, and it leads to the symplectic cohomology $SH^*(M)$ of
Cieliebak-Floer-Hofer \cite{cieliebak-floer-hofer95} and Viterbo
\cite{viterbo97a,viterbo97b}. Informally one can think of $SH^*(M)$
as the Floer cohomology $HF^*(M \setminus \partial M,H)$ for a
Hamiltonian function $H$ on the interior whose gradient points
outwards near the boundary, and becomes infinite as we approach the
boundary. For technical reasons, in the actual definition one takes
the direct limit over a class of functions with slower growth (to
clarify the conventions: our $SH^k(M)$ is dual to the $FH^{2n-k}(M)$
in \cite{viterbo97a}). The algebraic structure of symplectic
cohomology is different from the familiar case of closed $M$, where
one has large quantum cohomology and the WDVV equation. Operations
$SH^*(M)^{\otimes p} \rightarrow SH^*(M)^{\otimes q}$, for $p \geq 0$
and $q>0$, come from families of Riemann surfaces with $p+q$
punctures, together with a choice of local coordinate around each
puncture. The Riemann surfaces may degenerate to stable singular
ones, but only if no component of the normalization contains some of
the first $p$ and none of the last $q$ punctures. This means that if
we take only genus zero and $q = 1$ then no degenerations at all are
allowed, and the resulting structure is that of a Batalin-Vilkovisky
(BV) algebra \cite{getzler94}. For instance, let $M = D(T^*L)$ be a
unit cotangent bundle of an oriented closed manifold $L$. Viterbo
\cite{viterbo97b} computed that $SH^*(M) \iso H_{n-*}(\Lambda L)$ is
the homology of the free loop space, and a reasonable conjecture says
that the BV structure agrees with that of Chas-Sullivan
\cite{chas-sullivan99}.

Returning to the specific situation of Assumption
\ref{as:projective}, and supposing that $U$ has been chosen in such a
way that the Reeb flow on $\partial M$ becomes periodic, one can use
a Bott-Morse argument \cite{pozniak} to get a spectral sequence which
converges to $SH^*(M)$. The starting term is
\begin{equation} \label{eq:ss}
E_1^{pq} = \begin{cases}
 H^q(M) & p = 0, \\
 H^{q+3p}(\partial M) & p < 0.
\end{cases}
\end{equation}
It might be worth while to investigate this further, in order to
identify the differentials (very likely, a version of the relative
Gromov-Witten invariants \cite{ionel-parker98} for $D \subset X$).
But even without any more effort, one can conclude that each group
$SH^k(M)$ is finite-dimensional. In particular, assuming that
$dim_{\C}(X) > 2$ (and appealing to hard Lefschetz, which will be the
only time that we use any algebraic geometry) one has
\begin{equation} \label{eq:b2}
dim\, SH^2(M) \leq b_2(M) + b_0(\partial M) = b_2(X).
\end{equation}

\section{Fukaya categories}

$M$ (taken as in Assumption \ref{as:projective}) is an exact
symplectic manifold, and there is a well-defined notion of exact
Lagrangian submanifold in it. Such submanifolds $L$ have the property
that there are no non-constant holomorphic maps $u: (\Sigma,\partial
\Sigma) \rightarrow (M,L)$ for a compact Riemann surface $\Sigma$,
hence a theory of ``Gromov-Witten invariants with Lagrangian boundary
conditions'' would be trivial in this case. To get something
interesting, one removes some boundary points from $\Sigma$, thus
dividing the boundary into several components, and assigns different
$L$ to them. The part of this theory where $\Sigma$ is a disk gives
rise to the Fukaya $A_\infty$-category $\F(M)$.

The basic algebraic notion is as follows. An $A_\infty$-category $\A$
(over some field, let's say $\Q$) consists of a set of objects
$Ob\,\A$, and for any two objects a graded $\Q$-vector space of
morphisms $hom_\A(X_0,X_1)$, together with composition operations
\begin{align*}
 & \mu^1_\A : hom_\A(X_0,X_1) \longrightarrow hom_\A(X_0,X_1)[1], \\
 & \mu^2_\A : hom_\A(X_1,X_2) \otimes hom_\A(X_0,X_1) \longrightarrow
 hom_\A(X_0,X_2), \\
 & \mu^3_\A : hom_{\A}(X_2,X_3) \otimes hom_{\A}(X_1,X_2) \otimes
 hom_{\A}(X_0,X_1) \rightarrow \\ & \qquad \qquad \longrightarrow
 hom_{\A}(X_0,X_3)[-1], \quad \dots
\end{align*}
These must satisfy a sequence of quadratic ``associativity''
equations, which ensure that $\mu^1_\A$ is a differential, $\mu^2_\A$
a morphism of chain complexes, and so on. Note that by forgetting all
the $\mu^d_\A$ with $d \geq 3$ and passing to $\mu^1_\A$-cohomology
in degree zero, one obtains an ordinary $\Q$-linear category, the
induced cohomological category $H^0(\A)$ -- actually, in complete
generality $H^0(\A)$ may not have identity morphisms, but we will
always assume that this is the case (one says that $\A$ is
cohomologically unital).

\begin{figure}[htb] \begin{center}
\epsfig{file = 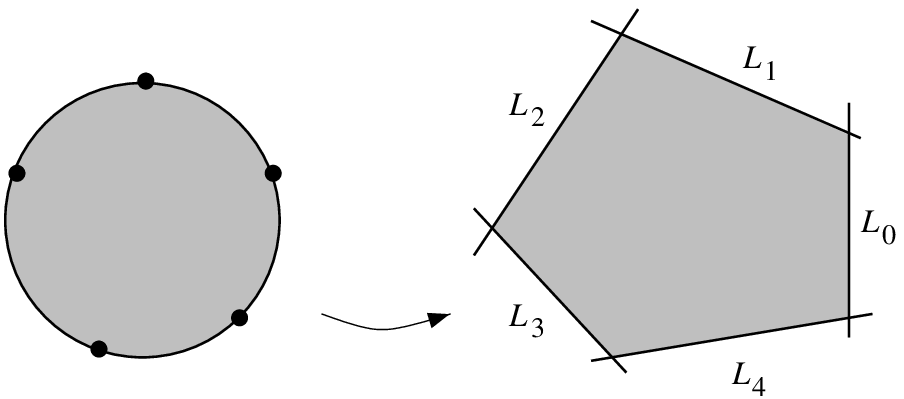}
\caption{\label{fig:polygon}}
\end{center} \end{figure}

In our application, objects of $\A = \F(M)$ are closed exact
Lagrangian submanifolds $L \subset M \setminus \partial M$, with a
bit of additional topological structure, namely a grading
\cite{kontsevich94,seidel99} and a $Spin$ structure \cite{fooo}. If
$L_0$ is transverse to $L_1$, the space of morphisms $hom_\A(L_0,L_1)
= CF(L_0,L_1)$ is generated by their intersection points, graded by
Maslov index. The composition $\mu_\A^d$ counts ``pseudo-holomorphic
$(d+1)$-gons'', which are holomorphic maps from the disk minus $d+1$
boundary points to $M$. The sides of the ``polygons'' lie on
Lagrangian submanifolds, and the corners are specified intersection
points; see Figure \ref{fig:polygon}. There are some technical issues
having to do with transversality, which can be solved by a small
inhomogeneous perturbation of the holomorphic map equation. This
works for all exact symplectic manifolds with contact type boundary,
satisfying $c_1 = 0$, and is quite an easy construction by today's
standards, since the exactness condition removes the most serious
problems (bubbling, obstructions).

It is worth while emphasizing that, unlike the case of Gromov-Witten
invariants, each one of the coefficients which make up $\mu^d_\A$
depends on the choice of perturbation. Only by looking at all of them
together does one get an object which is invariant up to a suitable
notion of quasi-isomorphism. To get something which is well-defined
in a strict sense, one can descend to the cohomological category
$H^0(\F(M))$ (which was considered by Donaldson before Fukaya's work)
whose morphisms are the Floer cohomology groups, with composition
given by the ``pair-of-pants'' product; but that is rather a waste of
information.

At this point, we must admit that there is essentially no chance of
computing $\F(M)$ explicitly. The reason is that we know too little
about exact Lagrangian submanifolds; indeed, this field contains some
of the hardest open questions in symplectic geometry. One way out of
this difficulty, proposed by Kontsevich \cite{kontsevich94}, is to
make the category more accessible by enlarging it, adding new objects
in a formal process, which resembles the introduction of chain
complexes over an additive category. This can be done for any
$A_\infty$-category $\A$, and the outcome is called the
$A_\infty$-category of twisted complexes, $Tw(\A)$. It contains the
original $A_\infty$-category as a full subcategory, but this
subcategory is not singled out intrinsically, and very different $\A$
can have the same $Tw(\A)$. The cohomological category $D^b(\A) =
H^0(Tw(\A))$, usually called the derived category of $\A$, is
triangulated (passage to cohomology is less damaging at this point,
since the triangulated structure allows one to recover many of the
higher order products on $Tw(\A)$ as Massey products). For our
purpose it is convenient to make another enlargement, which is
Karoubi or idempotent completion, and leads to a bigger
$A_\infty$-category $Tw^{\pi}(\A) \supset Tw(\A)$ and triangulated
category $D^\pi(\A) = H^0(Tw^{\pi}(\A))$. The main property of
$D^\pi(\A)$ is that for any object $X$ and idempotent endomorphism
$\pi: X \rightarrow X$, $\pi^2 = \pi$, there is a direct splitting $X
= im(\pi) \oplus ker(\pi)$. The details, which are not difficult,
will be explained elsewhere.

\section{Picard-Lefschetz theory}

We will now restrict the class of symplectic manifolds even further:

\begin{assumption} \label{as:lefschetz-fibre}
In the situation of Assumption \ref{as:projective}, suppose that $X$
is itself a hyperplane section in a smooth projective variety $Y$,
with $\K_Y \iso \O_Y(-X)$. Moreover, $X = X_0$ should be part of a
Lefschetz pencil of such sections $\{X_z\}$, whose base locus is $D =
X_0 \cap X_{\infty}$.
\end{assumption}

This gives a natural source of Lagrangian spheres in $M$, namely the
vanishing cycles of the Lefschetz pencil. Recall that to any
Lagrangian sphere $S$ one can associate a Dehn twist, or
Picard-Lefschetz monodromy map, which is a symplectic automorphism
$\tau_S$. The symplectic geometry of these maps is quite rich, and
contains information which is not visible on the topological level
\cite{seidel97,seidel98b,seidel99}. The action of $\tau_S$ on the
Fukaya category is encoded in an exact triangle in $Tw(\F(M))$, of
the form
\begin{equation} \label{eq:triangle}
\xymatrix{
 {L} \ar[rr] && {\tau_S(L)} \ar[dl]^-{[1]} \\
 & {\!\!\!\! HF^*(S,L) \otimes S \!\!\!\!} \ar[ul]
}
\end{equation}
for any $L$, and where the $\otimes$ is just a direct sum of several
copies of $S$ in various degrees. This is a consequence of the long
exact sequence in Floer cohomology \cite{seidel01}.

In the situation of Assumption \ref{as:lefschetz-fibre}, if we choose
a distinguished basis of vanishing cycles $S_1,\dots,S_m$ for the
pencil, the product of their Dehn twists is almost the identity map.
More precisely, taking into account the ``grading'' of the objects of
the Fukaya category, one finds that
\[
\tau_{S_1}\dots\tau_{S_m}(L) \iso L[2]
\]
where $[2]$ denotes change in the grading by 2. By combining this
trick with \eqref{eq:triangle} one can prove the following result:

\begin{thm} \label{th:generators}
$S_1,\dots,S_m$ are split-generators for $D^\pi(\F(M))$. This means
that any object of $Tw^\pi(\F(M))$ can be obtained from them, up to
quasi-isomorphism, by repeatedly forming mapping cones and idempotent
splittings.
\end{thm}

\section{Hochschild cohomology}

The Hochschild cohomology $HH^*(\A,\A)$ of an $A_\infty$-category
$\A$ can be defined by generalizing the Hochschild complex for
algebras in a straightforward way, or more elegantly using the
$A_\infty$-category $fun(\A,\A)$ of functors and natural
transformations, as endomorphisms of the identity functor. A
well-known rather imprecise principle says that ``Hochschild
cohomology is an invariant of the derived category''. In a rigorous
formulation which is suitable for our purpose,
\begin{equation} \label{eq:hh}
HH^*(\A,\A) \stackrel{?}{\iso} HH^*(Tw^\pi(\A),Tw^\pi(\A)).
\end{equation}
This is unproved at the moment, because $Tw^\pi(\A)$ itself has not
been considered in the literature before, but it seems highly
plausible (a closely related result has been proved in
\cite{keller98}). Hochschild cohomology is important for us because
of its role in deformation theory, see the next section; but we want
to discuss its possible geometric meaning first.

\begin{figure}[htb] \begin{center}
\epsfig{file=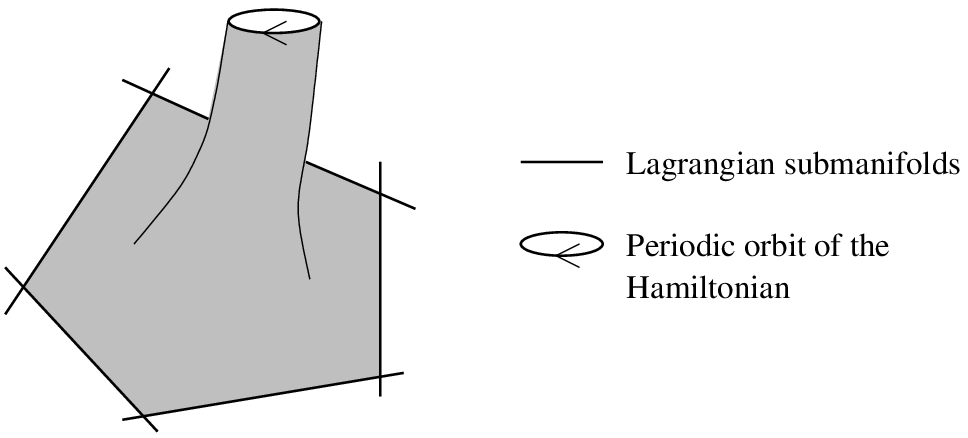} \caption{\label{fig:polygon2}}
\end{center} \end{figure}

Let $M$ be as in Assumption \ref{as:projective} (one could more
generally take any exact symplectic manifold with contact type
boundary and vanishing $c_1$). Then there is a natural ``open-closed
string map'' from the symplectic cohomology to the Hochschild
cohomology of the Fukaya category:
\begin{equation} \label{eq:open-closed}
SH^*(M) \longrightarrow HH^*(\F(M),\F(M)).
\end{equation}
This is defined in terms of Riemann surfaces obtained from the disk
by removing one interior point and an arbitrary number of boundary
points. Near the interior point, one deforms the holomorphic map
equation in the same way as in the definition of $SH^*(M)$, using a
large Hamiltonian function; otherwise, one uses boundary conditions
as for $\F(M)$. Figure \ref{fig:polygon2} shows what the solutions
look like.

$HH^*(\A,\A)$ for any $\A$ carries the structure of a Gerstenhaber
algebra, and one can verify that \eqref{eq:open-closed} is a morphism
of such algebras. Actually, since $SH^*(M)$ is a BV algebra, one
expects the same of $HH^*(\F(M),\F(M))$. This should follow from the
fact that $\F(M)$ is a cyclic $A_\infty$-category in some appropriate
weak sense, but the story has not yet been fully worked out (two
relevant papers for the algebraic side are \cite{tradler01} and
\cite{tamarkin-tsygan00}).

\begin{conjecture} \label{conj:hh}
If $M$ is as in Assumption \ref{as:lefschetz-fibre},
\eqref{eq:open-closed} is an isomorphism.
\end{conjecture}

Assumption \ref{as:lefschetz-fibre} appears here mainly for the sake
of caution. There are a number of cases which fall outside it, and to
which one would want to extend the conjecture, but it is not clear
where to draw the line. Certainly, without some restriction on the
geometry of $M$, there can be no connection between the Reeb flow on
$\partial M$ and Lagrangian submanifolds?

\section{Deformations of categories\label{sec:def}}

The following general definition, due to Kontsevich, satisfies the
need for a deformation theory of categories which should be
applicable to a wide range of situations: for instance, a deformation
of a complex manifold should induce a deformation of the associated
differential graded category of complexes of holomorphic vector
bundles. By thinking about this example, one quickly realizes that
such a notion of deformation must include a change in the set of
objects itself. The $A_\infty$-formalism, slightly extended in an
entirely natural way, fits that requirement perfectly. The relevance
to symplectic topology is less immediately obvious, but it plays a
central role in Fukaya, Oh, Ohta and Ono's work on ``obstructions''
in Floer cohomology \cite{fooo} (a good expository account from their
point of view is \cite{fukaya00}).

For concreteness we consider only $A_\infty$-deformations with one
formal parameter, that is to say over $\Q[[t]]$. Such a deformation
$\E$ is given by a set $Ob\,\E$ of objects, and for any two objects a
space $hom_\E(X_0,X_1)$ of morphisms which is a free graded
$\Q[[t]]$-module, together with composition operations as before but
now including a $0$-ary one: this consists of a so-called
``obstruction cocycle''
\begin{equation} \label{eq:obstruction}
\mu^0_\E \in hom_\E^2(X,X)
\end{equation}
for every object $X$, and it must be of order $t$ (no constant term).
There is a sequence of associativity equations, extending those of an
$A_\infty$-category by terms involving $\mu^0_\E$. Clearly, if one
sets $t = 0$ (by tensoring with $\Q$ over $\Q[[t]]$), $\mu^0_\E$
vanishes and the outcome is an ordinary $A_\infty$-category over
$\Q$. This is called the special fibre and denoted by $\E_{sp}$. One
says that $\E$ is a deformation of $\E_{sp}$.

A slightly more involved construction associates to $\E$ two other
$A_\infty$-categories, the global section category $\E_{gl}$ and the
generic fibre $\E_{gen}$, which are defined over $\Q[[t]]$ and over
the Laurent series ring $\Q[t^{-1}][[t]]$, respectively. One first
enlarges $\E$ to a bigger $A_\infty$-deformation $\E_c$ by coupling
the existing objects with formal connections (the terminology comes
from the application to complexes of vector bundles). Objects of
$\E_c$ are pairs $(X,\alpha)$ consisting of $X \in Ob\,\E$ and an
$\alpha \in hom^1_\E(X,X)$ which must be of order $t$. The morphism
spaces remain the same as in $\E$, but all the composition maps are
deformed by infinitely many contributions from the connection. For
instance,
\begin{equation} \label{eq:deformed-obstruction}
\mu^0_{\E_c} = \mu^0_\E + \mu^1_\E(\alpha) + \mu^2_\E(\alpha,\alpha)
+ \dots \in hom^2_{\E_c}((X,\alpha),(X,\alpha)) = hom^2_\E(X,X).
\end{equation}
$\E_{gl} \subset \E_c$ is the full $A_\infty$-subcategory of objects
for which \eqref{eq:deformed-obstruction} is zero; and $\E_{gen}$ is
obtained from this by inverting $t$. The transition from $\E_{sp}$ to
$\E_{gl}$ and $\E_{gen}$ affects the set of objects in the following
way: if for some $X$ one cannot find an $\alpha$ such that
\eqref{eq:deformed-obstruction} vanishes, then the object is
``obstructed'' and does not survive into $\E_{gl}$; if on the other
hand there are many different $\alpha$, a single $X$ can give rise to
a whole family of objects of $\E_{gl}$. Finally, two objects of
$\E_{gen}$ can be isomorphic even though the underlying objects of
$\E_{sp}$ aren't; this happens when the isomorphism involves negative
powers of $t$.

The classification of $A_\infty$-deformations of an
$A_\infty$-category $\A$ is governed by its Hochschild cohomology, or
rather by the dg Lie algebra underlying $HH^{*+1}(\A,\A)$, in the
sense of general deformation theory \cite{goldman-millson88}. We
cannot summarize that theory here, but as a simple example, suppose
that $HH^2(\A,\A) \iso \Q$. Then a nontrivial $A_\infty$-deformation
of $\A$, if it exists, is unique up to equivalence and change of
parameter $t \mapsto f(t)$ (to be accurate, $f(t)$ may contain roots
of $t$, so the statement holds over $\Q[[t,t^{1/2},t^{1/3},\dots]]$).
The intuitive picture is that the ``versal deformation space'' has
dimension $\leq 1$, so that any two non-constant arcs in it must
agree up to reparametrization.

In the situation of Assumption \ref{as:projective}, the embedding of
our exact symplectic manifold $M$ into $X$ should give rise to an
$A_\infty$-deformation $\F(M \subset X)$. We say ``should'' because
the details, which in general require the techniques of \cite{fooo},
have not been carried out yet. Roughly speaking one takes the same
objects as in $\F(M)$ and the same morphism spaces, tensored with
$\Q[[t]]$, but now one allows ``holomorphic polygons'' which map to
$X$, hence may intersect the divisor $D$. The numbers of such
polygons intersecting $D$ with multiplicity $k$ will form the $t^k$
term of the composition maps in $\F(M \subset X)$. Because there can
be holomorphic discs bounding our Lagrangian submanifolds in $X$,
nontrivial obstruction cocycles \eqref{eq:obstruction} may appear.

The intended role of $\F(M \subset X)$ is to interpolate between
$\F(M)$, which we have been mostly discussing up to now, and the
Fukaya category $\F(X)$ of the closed symplectic manifold $X$ as
defined in \cite{fooo,fukaya01b}. The $t^0$ coefficients count
polygons which are disjoint from $D$, and these will automatically
lie in $M$, so that
\[
\F(M \subset X)_{sp} \iso \F(M).
\]
The relation between the generic fibre and $\F(X)$ is less
straightforward. First of all, $\F(M \subset X)_{gen}$ will be an
$A_\infty$-category over $\Q[t^{-1},t]]$, whereas $\F(X)$ is defined
over the Novikov ring $\Lambda_t$. Intuitively, one can think of this
difference as the consequence of a singular deformation of the
symplectic form. Namely, if one takes a sequence of symplectic forms
(all in the same cohomology class) converging towards the current
$[D]$, the symplectic areas of holomorphic discs $u$ would tend to
the intersection number $u \cdot D$. A more serious issue is that
$\F(M \subset X)_{gen}$ is clearly smaller than $\F(X)$, because it
contains only Lagrangian submanifolds which lie in $M$. However, that
difference may disappear if one passes to derived categories:

\begin{conjecture} \label{conj:generic-fibre} In the situation of
Assumption \ref{as:lefschetz-fibre}, there is a canonical equivalence
of triangulated categories
\[
D^\pi(\F(M \subset X)_{gen} \otimes_{\Q[t^{-1}][[t]]} \Lambda_t) \iso
D^\pi(\F(X)).
\]
\end{conjecture}

In comparison with the previous conjecture, Assumption
\ref{as:lefschetz-fibre} is far more important here. The idea is that
there should be an analogue of Theorem \ref{th:generators} for
$D^\pi(\F(X))$, saying that this category is split-generated by
vanishing cycles, hence by objects which are also present in $\F(M
\subset X)$.

To pull together the various speculations, suppose that $Y =
\CP{n+1}$ for some $n \geq 3$; $X \subset Y$ is a hypersurface of
degree $n+2$; and $D \subset X$ is the intersection of two such
hypersurfaces. Then $D^\pi(\F(M))$ is split-generated by finitely
many objects, hence $Tw^\pi(\F(M))$ is at least in principle
accessible to computation. Conjecture \ref{conj:hh} together with
\eqref{eq:b2}, \eqref{eq:hh} tells us that $HH^2(\F(M),\F(M)) \iso
HH^2(Tw^\pi(\F(M)),Tw^\pi(\F(M)))$ is at most one-dimensional, so an
$A_\infty$-deformation of $Tw^\pi(\F(M))$ is unique up to a change of
the parameter $t$. From this deformation, Conjecture
\ref{conj:generic-fibre} would enable one to find $D^\pi(\F(X))$,
again with the indeterminacy in the parameter (fixing this is
somewhat like computing the mirror map).
\providecommand{\bysame}{\leavevmode\hbox
to3em{\hrulefill}\thinspace}
\providecommand{\MR}{\relax\ifhmode\unskip\space\fi MR }
\providecommand{\MRhref}[2]{%
  \href{http://www.ams.org/mathscinet-getitem?mr=#1}{#2}
} \providecommand{\href}[2]{#2}


\begin{thebibliography}{10}

\bibitem{chas-sullivan99}
M.~Chas and D.~Sullivan, \emph{String topology}, Preprint
math.GT/9911159.

\bibitem{cieliebak-floer-hofer95}
K.~Cieliebak, A.~Floer, and H.~Hofer, \emph{Symplectic homology {II}:
a general
  construction}, Math. Z. \textbf{218} (1995), 103--122.

\bibitem{eliashberg-givental-hofer02}
Ya. Eliashberg, A.~Givental, and H.~Hofer, \emph{Introduction to
symplectic
  field theory}, Geom. Funct. Anal. \textbf{Special Volume, Part II} (2000),
  560--673.

\bibitem{fukaya02}
K.~Fukaya, \emph{Asymptotic analysis, multivalued {M}orse theory, and
mirror
  symmetry}, Preprint 2002.

\bibitem{fukaya00}
\bysame, \emph{Deformation theory, homological algebra, and mirror
symmetry},
  Preprint, December 2001.

\bibitem{fukaya01b}
\bysame, \emph{Floer homology and mirror symmetry {II}}, Preprint
2001.

\bibitem{fukaya93}
\bysame, \emph{Morse homotopy, {$A_\infty$}-categories, and {F}loer
  homologies}, Proceedings of {GARC} workshop on Geometry and Topology (H.~J.
  Kim, ed.), Seoul National University, 1993.

\bibitem{fukaya-oh98}
K.~Fukaya and Y.-G. Oh, \emph{Zero-loop open strings in the cotangent
bundle
  and {M}orse homotopy}, Asian J. Math. \textbf{1} (1998), 96--180.

\bibitem{fooo}
K.~Fukaya, Y.-G. Oh, H.~Ohta, and K.~Ono, \emph{Lagrangian
intersection {F}loer
  theory - anomaly and obstruction}, Preprint, 2000.

\bibitem{getzler94}
E.~Getzler, \emph{{B}atalin-{V}ilkovisky algebras and 2d
{T}opological {F}ield
  {T}heories}, Commun. Math. Phys \textbf{159} (1994), 265--285.

\bibitem{goldman-millson88}
W.~Goldman and J.~Millson, \emph{The deformation theory of the
fundamental
  group of compact {K}{\"a}hler manifolds}, IHES Publ. Math. \textbf{67},
  43--96.

\bibitem{ionel-parker98}
E.-N.~Ionel and T.~Parker, \emph{Gromov-{W}itten invariants of
symplectic
  sums}, Math. Res. Lett. \textbf{5} (1998), 563--576.

\bibitem{keller98}
B.~Keller, \emph{Invariance and localization for cyclic homology of
{DG}
  algebras}, J. Pure Appl. Alg. \textbf{123} (1998), 223--273.

\bibitem{kontsevich94}
M.~Kontsevich, \emph{Homological algebra of mirror symmetry},
Proceedings of
  the International Congress of Mathematicians (Z{\"u}rich, 1994),
  Birkh{\"a}user, 1995, pp.~120--139.

\bibitem{kontsevich-soibelman00}
M.~Kontsevich and Y.~Soibelman, \emph{Homological mirror symmetry and
torus
  fibrations}, Symplectic geometry and mirror symmetry, World Scientific, 2001,
  pp.~203--263.

\bibitem{li-ruan98}
A.-M. Li and Y.~Ruan, \emph{Symplectic surgery and {G}romov-{W}itten
invariants
  of {C}alabi-{Y}au 3-folds}, Invent. Math. \textbf{145} (2001), 151--218.

\bibitem{ozsvath-szabo01}
P.~Ozsvath and Z.~Szabo, \emph{Holomorphic disks and topological
invariants for
  rational homology three-spheres}, Preprint math.SG/0101206.

\bibitem{polishchuk-zaslow98}
A.~Polishchuk and E.~Zaslow, \emph{Categorical mirror symmetry: the
elliptic
  curve}, Adv. Theor. Math. Phys. \textbf{2} (1998), 443--470.

\bibitem{pozniak}
M.~Po{\'z}niak, \emph{Floer homology, {N}ovikov rings and clean
intersections},
  Northern California Symplectic Geometry Seminar, Amer. Math. Soc., 1999,
  pp.~119--181.

\bibitem{seidel97}
P.~Seidel, \emph{Floer homology and the symplectic isotopy problem},
Ph.D.
  thesis, Oxford University, 1997.

\bibitem{seidel98b}
\bysame, \emph{Lagrangian two-spheres can be symplectically knotted},
J.
  Differential Geom. \textbf{52} (1999), 145--171.

\bibitem{seidel99}
\bysame, \emph{Graded {L}agrangian submanifolds}, Bull. Soc. Math.
France
  \textbf{128} (2000), 103--146.

\bibitem{seidel01}
\bysame, \emph{A long exact sequence for symplectic {F}loer
cohomology},
  Preprint math.SG/0105186.

\bibitem{tamarkin-tsygan00}
D.~Tamarkin and B.~Tsygan, \emph{Noncommutative differential
calculus, homotopy
  {BV} algebras and formality conjectures}, Preprint math.KT/0002116.

\bibitem{tradler01}
T.~Tradler, \emph{Infinity-inner-products on {A}-infinity-algebras},
Preprint
  math.\-AT/0108027.

\bibitem{viterbo97a}
C.~Viterbo, \emph{Functors and computations in {F}loer homology with
applications,
  {P}art {I}}, Geom. Funct. Anal. \textbf{9} (1999), 985--1033.

\bibitem{viterbo97b}
\bysame, \emph{Functors and computations in {F}loer homology with
  applications, {P}art {II}}, Preprint 1996.
\end{thebibliography}
\end{document}